# Investigations of the limit distribution and the asymptotic behavior of summation arithmetic functions

VICTOR VOLFSON

ABSTRACT. The aim of the paper is to study the limit distributions and the asymptotic behavior of summation arithmetic functions. A probabilistic approach based on the use of the axioms of probability theory is used for these purposes. Sufficient conditions are proved under which these functions have a limiting normal distribution. Arithmetic functions having a limiting normal distribution are found in the paper. The author investigates summation functions of a general form and finds sufficient conditions under which they have a limiting normal distribution. Examples of arithmetic functions that satisfy these requirements are considered. We prove an ergodic theorem for summation arithmetic functions, for which the sequence of random variables has the stationarity property in the broad sense. The asymptotics of the growth of the deviation of the values of the arithmetic function from its mean value are investigated. It is shown that an equivalent formulation of the Riemann hypothesis for the Mertens function is satisfied almost everywhere.

1. INTRODUCTION

An arithmetic function in the general case is a function defined on the set of natural numbers and taking values on the set of complex numbers. The name arithmetic function is related to the fact that this function expresses some arithmetic property of the natural series.

Many arithmetic functions, considered in number theory, take on natural values. Such functions include functions of the number of natural numbers:

$$Q(n) = \sum_{k \in S} 1, \qquad (1.1)$$

satisfying a certain property.

---

An arithmetic function of an adder type is a function:

$$S(x) = \sum_{n \leq x} f(n). \tag{1.2}$$

A special case of an arithmetic function of an adder type is the functions (1.1).

We usually limited ourselves to the average value of the arithmetic function on the interval of the natural series $[1, n]$:

$$M_f(n) = \frac{1}{n} \sum_{k=1}^{n} f(k) \tag{1.3}$$

when studying the distribution of the values of arithmetic functions in classical and even more recent studies and looked his an asymptotic approximate expression in terms of the simplest possible arithmetic functions.

However, it is obvious that the values of the arithmetic function can oscillate about the average value (1.3) over very wide limits even in the case of simple arithmetic functions. Therefore, the natural question is how much the value of the arithmetic function deviates from its mean value for large values of the argument. It is possible to pose a broader question about the distribution of an arithmetic function on a natural number for large values, and hence, what is the limiting distribution of the arithmetic function.

A powerful method for solving problems of the distribution of the values of a function of a natural argument is the method of trigonometric sums. Vinogradov carried out the first statements of problems on the distribution of values of arithmetic functions of the summation type (1.2) and their solution using the trigonometric sum method [1].This line of research was continued [2]. However, in general, the method of trigonometric sums is used to find the distribution of arithmetic functions possessing the additivity or multiplicative property.

The axiomatic of probability theory makes it possible to use probability theory to find the distribution of any arithmetic functions. Any interval of a natural number $[1, n]$ can be naturally transformed into a probability space - $(\Omega_n, A_n, P_n)$, taking for $\Omega_n = \{1, 2, ..., n\}$, $A_n$ - all subsets $\Omega_n$, $P_n(A) = N(k \in A)/n$, where $N(k \in A)$ is the number of natural numbers in a subset $A$.

Then an arbitrary real arithmetic function $f(k)$ on $\Omega_n$ can be considered as a random variable $x_n$ on this probability space:

$$x_n(k) = f(k)(1 \leq k \leq n). \tag{1.4}$$

One can speak, in particular, of the mean value of the arithmetic function (mathematical expectation) $M(x_n) = M_f(n)$, determined by formula (1.3), variance:

$$D(x_n) = D_f(n) = \frac{1}{n}\sum_{k=1}^{n} f^2(k) - (\frac{1}{n}\sum_{k=1}^{n} f(k))^2$$

and the distribution function:

$$F_n(y) = P_n\{f(k) < y\}.$$

Thus, we can associate for each real arithmetic function $f(n)$ a sequence of random variables $x_n$ that are on different probability spaces.

This probabilistic approach is used to the question of the distribution of arithmetic functions: [3], [4]. However, these papers are mainly devoted to the consideration of additive and multiplicative arithmetic functions. Therefore, a large number of arithmetic functions remain unexplored.

The aim of this paper is to study the limit distributions and the asymptotic behavior of summation arithmetic functions. A probabilistic approach will be used for these purposes, based on axioms of the probability theory.

We call deterministic (defined) arithmetic functions, for which the deviation from the mean value is zero. Accordingly, arithmetic functions are called non-deterministic, for which the deviation from the mean value is not zero. We will consider non-deterministic functions in the paper.

Examples of summative nondeterministic arithmetic functions are: the number of prime numbers, the number of prime tuples in the natural series, the number of natural numbers for which the polynomials take prime values, the number of natural numbers having prime divisors, the number of natural numbers free of squares, the number of natural numbers, having an even or odd number of prime divisors, the Mertens function, Chebyshev functions, etc.

## 2. PROPERTIES OF A SEQUENCE OF RANDOM VAIABLES

We investigate summation arithmetic functions, which can be represented as sums of random variables $x_i$:

$$S_n = \sum_{i=1}^{n} x_i. \tag{2.1}$$

We recall some properties of a sequence of random variables $x_i (i = 1,...,n)$. The sequence of random Variables $x_i (i = 1,...,n)$ in (2.1) is stationary in the wide sense if the following conditions are satisfied:

1. $M(x_n)$ does not depend on $n$. (2.2)

2. $\mathrm{cov}(x_{n+m}, x_n) = M((x_{n+m} - M(x_{n+m}))(x_n - M(x_n)))$ does not depend on $n$. (2.3)

3. $D(x_n) < \infty$. (2.4)

We are interested in the asymptotic behavior of arithmetic functions, therefore conditions (2.2), (2.3), (2.4) must be satisfied for $n \to \infty$.

Based on (1.4) $x_n(k) = f(k)(1 \leq k \leq n)$, therefore:

$$\lim_{n \to \infty} \frac{S(n)}{n} = \lim_{n \to \infty} \frac{\sum_{k=1}^{n} f(k)}{n} = \lim_{n \to \infty} \frac{\sum_{k=1}^{n} x_n(k)}{n} = \lim_{n \to \infty} M(x_n) \tag{2.5}$$

If the condition (2.2) is satisfied asymptotically, having in mind (2.5), the following equality must hold:

$$\lim_{n \to \infty} M(x_n) = \lim_{n \to \infty} \frac{S(n)}{n} = C, \tag{2.6}$$

where $C$ is a constant.

Now let us verify the fulfillment of condition (2.4):

$$D(x_n) = M(x_n^2) - M^2(x_n) \leq M(x_n^2) = \frac{\sum_{k=1}^{n} x_n^2(k)}{n} \leq \sup_k |x_n^2(k)|. \tag{2.7}$$

Therefore, if the value $x_n$ is bounded, then (2.4) is satisfied on the basis of (2.7).

The sequence of random variables $x_1, x_2, ...$ is independent if:

$$P(x_1 \in B_1, x_2 \in B_2,...) = P(x_1 \in B_1)P(x_2 \in B_2)... \qquad (2.8)$$

The random variables $x_k$ and $x_{k+n}$ are asymptotically independent (for $n \to \infty$), when for any values $B_1, B_2 \in B$:

$$P(x_k \in B_1, x_{k+n} \in B_2) \to P(x_0 \in B_1)P(x_0 \in B_2). \qquad (2.9)$$

If the random variables of the sequence $x_i (i = 1,...,n)$ are asymptotically independent, then:

$$\lim_{n \to \infty} \text{cov}(x_k, x_{k+n}) = 0, \qquad (2.10)$$

(2.10) does not depend on n and (2.3) is satisfied.

Now let's talk about the property of strong mixing for stationary sequences in a broad sense.

Let's start with the definition. A stationary sequence in a broad sense $x_n$ satisfies the condition of strong mixing if:

$$\alpha(n) = \sup\nolimits_{A \in M_1^n, B \in M_{n+k}^{\infty}} |P(AB) - P(A)P(B)| \to 0 \text{ with } n \to \infty, \qquad (2.11)$$

where $\alpha(n)$ is the mixing coefficient, $A \in M_1^n$ means that the event $A$ belongs to the algebra of events $M_{n+k}^{\infty}$, similarly with respect to the event $B$. Later we consider several theorems based on the condition of strong mixing.

Let us examine the summation arithmetic functions from the point of view of the properties of sequences of random variables indicated above $x_1, x_2,...$. Consideration begins with an arithmetic function of the number of natural numbers $Q(n)$ that satisfy a certain property.

A random variable $x_k$ in a sequence $x_1, x_2,...$ for a function $Q(n)$ takes only two values: $x_k = 1$ if the natural number $k$ satisfies a certain property and $x_k = 0$ - otherwise. Thus, the random variables $x_k$ are bounded in this case; therefore, on based on (2.7), the condition (2.4) of stationarity in the broad sense is satisfied for functions $Q(n)$.

Sequences of random variables $x_1, x_2, \ldots$ for Chebyshev arithmetic functions also take only two values: 0 and $\ln(p)$. However, in connection with the fact that the second value is unlimited, then the stationarity condition in the broad sense is not satisfied for these arithmetic functions.

Now we check condition (2.6) for arithmetic functions $Q(n)$.

Based on [6], the formula for the arithmetic function of the number of square-free numbers is hold:

$$Q(n) = 6n/\pi^2 + O(n^{1/2}).$$

Therefore, this condition is fulfilled:

$$\lim_{n \to \infty} \frac{Q(n)}{n} = 6/\pi^2.$$

Let us check the condition for the asymptotic independence of the corresponding sequences for arithmetic functions $Q(n)$ - the last condition of stationarity in the broad sense.

Let's look at the arithmetic function of the number of primes. Suppose, if the natural $k (k > 2)$ is a prime, then the value of the random variable $x_k = 1$, otherwise - $x_k = 0$. Similarly with respect to $x_{k+1}$. We denote: $P(x_k = 1) = a (0 \leq a \leq 1)$, $P(x_k = 0) = 1 - a$, $P(x_{k+1} = 1) = b (0 \leq b \leq 1)$, $P(x_{k+1} = 0) = 1 - b$. Then, based on (2.8) and (2.9), on the one hand:

$P(x_k = 1, x_{k+1} = 0) = P(x_k = 1)P(x_{k+1} = 0 / x_k = 1) = a \cdot 1 = a$, because if $x_k = 1$, then $k$ - an odd number. On the other hand: $P(x_k = 1)P(x_{k+1} = 0) = a(1 - b)$. Therefore: $P(x_k = 1, x_{k+1} = 0) = a$ is not equal to $P(x_k = 1)P(x_{k+1} = 0) = a(1 - b)$, if $b$ not equal 0.

Thus, random variables $x_k, x_{k+1}$ are dependent. This situation is preserved when $n \to \infty$. Thus, in this case the asymptotic independence condition is not satisfied. If we take the arithmetic functions of the number of natural numbers having two (small) number of prime divisors, then the indicated dependence on the fact that the natural number $k$ is an even or odd number, as for prime numbers, is preserved.

There is an independence from the fact that the natural number $k$ is an even or odd number for arithmetic functions of the number of natural numbers having a large number of prime divisors (for large values $n$), therefore asymptotic independence holds.

Also asymptotic independence is performed for summation arithmetic functions when the defining property of which is the presence of an even or odd number of prime divisors for a natural number. Such arithmetic functions include: the number of natural numbers having simple divisors, where is a large number, the number of natural numbers that are free from squares, the number of natural numbers having an even or odd number of prime divisors, the Martens function, and others.

The stationarity in the broad sense is satisfied for all considered arithmetic functions $Q(n)$. The exceptions are arithmetic functions of the number of primes and the number of natural numbers with a small number of prime divisors.

The stationarity in the broad sense is satisfied in the case when the functions satisfy conditions (2.2), (2.3) for bounded summation arithmetic functions other than those indicated above and taking more than two discrete values.

For example, it was said above that the Mertens function satisfies the condition of asymptotic independence of the corresponding random variables; for this function the condition (2.3) is satisfied. It is known that for the Mertens function the following relation holds $M(n) = o(n)$, and therefore (2.6) holds:

$$\lim_{n \to \infty} \frac{M(n)}{n} = 0 \qquad (2.12)$$

and hence condition (2.3). Thus, the Mertens function satisfies all the stationary conditions in the broad sense.

3. THE LIMITING DISTRIBUTION FOR THE ARITHMETIC FUNCTIOM OF THE NUMBER OF NATURAL NUMBERS $Q(n)$

Let us consider the arithmetic function of the number of natural numbers having a certain property - $Q(n)$. The arithmetic function is the mapping: $Q(n): \mathbb{N} \to \mathbb{N}$.

$Q(n)$ is a piecewise constant, "stepwise", monotonically increasing, unbounded function. On the other hand the value is $0 \leq Q(n) \leq n$. The value $Q(n)$ increases by 1 if the natural satisfies a certain property. Otherwise, the value of the function does not change.

Examples of $Q(n)$ are arithmetic functions: the number of square-free natural numbers; number of prime numbers, etc.

Let us formulate the problem of finding the limiting distribution for an arithmetic function $Q(n)$.

We introduce probabilistic spaces: $(\Omega_n, \mathcal{A}_n, \mathbb{P}_n)$ by taking $\Omega_n = \{1, 2, ..., n\}$, $\mathcal{A}_n$ - all subsets $\Omega_n$, $P_n(A) = \frac{1}{n}\{N(m \in A\}$, where $N(m \in A)$ is the number of terms of the natural series satisfying the condition $m \in A$.

Then the arithmetic function $Q(m)$ can be considered as a random variable $S_n(m) = Q(m)(1 \leq m \leq n)$ on the indicated spaces.

We fix $n$ and consider the probability space of a random variable $S_n(m) = Q(m)$.

A computational arithmetic function $Q(n)$ can be represented as a sum of random variables on a given probability space:

$$S_n = \sum_{i=1}^{n} x_i, \qquad (3.1)$$

where $x_i = 1$ if a natural number $i$ has a certain property and $x_i = 0$ if the natural number does not have a certain property.

The probability of the event that $x_i = 1$ in (3.1) in a given probability space is determined by the formula:

$$p_n = \frac{1}{n} N\{x_i = 1\}. \qquad (3.2)$$

The probability that $x_i = 0$ is equal to $1 - p_n$.

Having in mind $N\{x_i = 1\} = Q(n)$, therefore (3.2) can be written in the form:

$$p_n = \frac{1}{n} N\{x_i = 1\} = Q(n)/n. \qquad (3.3)$$

First we assume that the random variables $x_i$ are independent.

Then the random variable $S_n$ has a binomial distribution. The limiting distribution for a binomial distribution is the Poisson distribution, or the normal distribution, depending on the behavior of the variance at the value $n \to \infty$:

$$D(S_n) = np_n(1 - p_n) = Q(n)(1 - Q(n)/n). \qquad (3.4)$$

Taking into account that $Q(n)$ increases without limit for the value $n \to \infty$ and $0 \leq 1 - Q(n)/n \leq 1$, then based on (3.4) the value $D(S_n)$ increases indefinitely. Therefore, based on the integral theorem of Moivre-Laplace, the limiting distribution for the arithmetic function $Q(n)$ is the normal distribution in this case.

Hardy-Littlewood's hypotheses about prime tuples and Betman-Horn consider that events corresponding to the appearance of a tuple of primes in the natural number and the natural number at which the polynomials in the Beitman-Horn hypothesis take on prime meanings are independent. Therefore, based on the assertion proved above, the arithmetic functions of the number of prime tuples in the natural series and the number of natural numbers for which the polynomials take prime values have a limiting normal distribution. It is possible to give other examples of arithmetic functions $Q(n)$ that are the sum of independent random variables, and therefore, based on the previously proved statement having a limiting normal distribution.

Now we consider the case when the arithmetic function $Q(n)$ is the sum of dependent random variables. We pose the question - what sufficient conditions must be satisfied in this case so that the arithmetic function $Q(n)$ has a limiting normal distribution.

The theorem is proved in [5]. Suppose that a sequence of random variables $x_n (n = 1, 2, ...)$ stationary in the wide sense satisfies the conditions of strong mixing with $\sum_{n=1}^{\infty} \alpha(n) < \infty$. Suppose, in addition, that the random variable is bounded with probability 1 - $|x_n| = c_o < \infty$. Then, if $D(x_n) \neq 0$, then the random variable $S_n$ has a normal distribution.

Arithmetic functions were introduced in Chapter 2, for which the stationarity conditions in the broad sense are satisfied for the corresponding sequence of random variables $x_n (n=1,2,...)$. In addition, for a random variable $x_n$ corresponding to $Q(n)$, the boundedness condition is satisfied $|x_n| \leq 1$.

Based on (3.3), the formula for the variance of the random variables $x_n$ (corresponding to $Q(n)$) is true

$$D(x_n) = p_n(1-p_n) = \frac{Q(n)}{n}(1-\frac{Q(n)}{n}) \neq 0,$$

if $Q(n) \neq n$.

Consequently, it remains to verify the strong mixing property with $\sum_{n=1}^{\infty} \alpha(n) < \infty$ for sequences $x_n (n=1,2,...)$ corresponding to the indicated arithmetic functions $Q(n)$ in order to verify the sufficient conditions for the limit theorem.

As examples, we consider arithmetic functions $Q(n)$: the number of prime numbers and the number of natural numbers that are free from squares.

Based on (2.11), we find the strong mixing coefficient for a random sequence corresponding to the arithmetic function of the number of prime numbers. Let us take the event $A$ that the number $k$ is prime, i.e. $x_k = 1$. Let us take the event $B$ that the number $n+k$ is composite, i.e. $x_{n+k} = 0$. Then $P(x_k=1) = 1/\ln k + o(1/\ln k)$, $P(x_{n+k}=0) = 1 - 1/\ln(n+k) + o(1/\ln(n+k))$ and $P(A)P(B) = (1 + 1/\ln k + o(1/\ln k))(1 - 1/\ln(n+k) + o(1/\ln(n+k)))$. On the other hand $P(AB) = P(A)P(B/A) = 1/\ln k + o(1/\ln k)$, since $P(B/A) = 1$, if the number $n+k$ is even. The last assumption is made to obtain the upper bound of the difference. Therefore, the strong mixing ratio in this case is:

$$\alpha(n) = \sup |P(AB) - P(A)P(B)| = (1/\ln k + o(1/\ln k))(1/\ln(n+k) + o(1/\ln(n+k))). \quad (3.5)$$

Having in mind (3.5) $\alpha(n) \to 0$ with $n \to \infty$, therefore the strong mixing condition is satisfied in this case.

However, it is required in the above theorem that $\sum_{n=1}^{\infty} \alpha(n) < \infty$. Let's check the fulfillment of this condition for the arithmetic function of the number of primes:

$$\sum_{n=1}^{\infty} \alpha(n) \geq [1/\ln k + o(1/\ln k)] \sum_{n=1}^{\infty} 1/\ln(n+k) = \infty. \tag{3.6}$$

Based on (3.6), this sufficient condition for the limit theorem for the arithmetic function of the number of primes is not satisfied. In addition, the stationary condition in the broad sense is also not satisfied for this function.

Now we find the strong mixing coefficient for a sequence of random variables of the corresponding arithmetic function of the number of natural numbers free from squares. Based on [6], the following formula holds for the indicated arithmetic function:

$$Q(n) = 6n/\pi^2 + O(n^{1/2}).$$

We take the event $A$ that the number $k$ is free from squares, i.e. $x_k = 1$. We take the event $B$ that the number $n+k$ is not free of squares, i.e. $x_{n+k} = 0$. Then $P(x_k = 1) = 6/\pi^2 + O(1/k^{1/2})$, $P(x_{n+k} = 0) = 1 - 6/\pi^2 + O(1/(n+k)^{1/2})$. The product of probabilities is equal to: $P(A)P(B) = [6/\pi^2 + O(1/k^{1/2})][1 - 6/\pi^2 + O(1/(n+k)^{1/2})]$ in this case. Then the value $P(AB)$ is: $P(A)P(B/A) = [6/\pi^2 + O(1/k^{1/2})][1 - 6/\pi^2 + O(1/(n+k)^{1/2})]$. The last equality is explained by the independence of the event $B$ from the event $A$. Therefore, the strong mixing condition is satisfied in this case:

$$\alpha(n) = \sup | P(AB) - P(A)P(B) | = 0. \tag{3.7}$$

Based on (3.7): $\sum_{n=1}^{\infty} \alpha(n) = 0 < \infty$, therefore, all the conditions of the above limit theorem are satisfied for a sequence of random variables corresponding to a given arithmetic function and the limiting distribution for a given function is normal. This is consistent with the QEIS A158819 schedule.

The conditions of the above limit theorem also hold for sequences of random variables corresponding to the following arithmetic functions $Q(n)$: the number of natural numbers having $k$ prime divisors, where $k$ is a large number, the number of natural numbers having an even or odd number of prime divisors. Therefore, the limiting distribution for them is also normal.

## 4. A LIMIT DTSTRIBUTION FOR SUMMATION ARITHMETIC FUNCTIONS IN GENERAL CASE

Any arithmetic summation function has the form:

$$S(n) = \sum_{i=1}^{n} f(i), \qquad (4.1)$$

where $f(n)$ is the arithmetic function. Based on (4.1), the summation arithmetic function is a piecewise constant, "step", unbounded function with a value $n \to \infty$. The step height depends on the value $f(i)$ at the point $i$.

Examples of summation functions are already considered in Chapter 3 arithmetic functions of the number of natural numbers that satisfy a certain property - $Q(n)$, the Mertens function, Chebyshev functions and others.

We pose the problem of finding the limit distribution for the summation function $S(n)$ in the general case.

We introduce probabilistic spaces $(\Omega_n, \mathcal{A}_n, \mathbb{P}_n)$: by taking $\Omega_n = \{1, 2, ..., n\}$, $\mathcal{A}_n$ - all subsets $\Omega_n$, $P_n(A) = \frac{1}{n}\{N(m \in A)\}$, where $N(m \in A)$ is the number of terms of the natural series satisfying the condition $m \in A$.

Then the arithmetic function $S(m)$ can be considered as a random variable $S_n(m) = S(m)(1 \leq m \leq n)$ on the indicated spaces.

We fix $n$ and consider the probability space of a random variable $S_n(m) = S(m)$.

A summation arithmetic function $S(n)$ can be represented as a sum of random variables on a given probability space:

$$S_n = \sum_{i=1}^{n} x_i, \qquad (4.2)$$

where a random variable $x_i$ can take, in general, more than two different discrete values.

Earlier we defined the class of arithmetic functions for which the stationary properties in the broad sense are satisfied. We also gave a theorem, under the conditions of which, the limiting distribution of an arithmetic function is a normal distribution. We rephrase it a little for the summation arithmetic function of the general form.

If the sequence of random variables $x_n (n=1,2,...)$ corresponding to the summation arithmetic function $S(n)$ is stationary in the broad sense and it has the strong mixing property with $\sum_{n=1}^{\infty} \alpha(n) < \infty$, $x_n$ is bounded with probability 1 and $D(x_n) \neq 0$, then the summation arithmetic function has a limiting normal distribution. Thus, we must verify these conditions for the summation arithmetic function under study.

Let's consider, as an example, the arithmetic summation function $S(n) = \sum_{i=1}^{n} f(i)$, where, $f(i) = 2$ if $i$ has an even number of prime divisors, $f(i) = 0$ if $i$ is not a square-free number and $f(i) = -1$ if $i$ has an odd number of prime divisors.

Let us first verify the stationarity condition in the broad sense of the corresponding sequence of random discrete quantities. It is known [7] that the probability of a natural number $n$ to have an even number of prime divisors is equal to the probability of a natural number $n$ to have an odd number of prime divisors:

$$v_1(n) = v_3(n) = 3/\pi^2 + o(1). \tag{4.3}$$

Therefore, the probability of a natural number $n$ hasn't prime divisors that are free from squares is equal to:

$$v_2(n) = 1 - 6/\pi^2 + o(1). \tag{4.4}$$

Having in mind (4.3) and (4.4), the corresponding limiting probabilities are equal to:

$$\begin{aligned} \lim_{n\to\infty} v_1(n) &= 3/\pi^2 \\ \lim_{n\to\infty} v_2(n) &= 1 - 6/\pi^2 \\ \lim_{n\to\infty} v_3(n) &= 3/\pi^2 \end{aligned} \tag{4.5}$$

Based on (4.3) - (4.5), the mathematical expectation of a random variable:

$$M(x_n) = 3/\pi^2 + o(1) \ . \tag{4.6}$$

Having in mind (4.6) the following holds:

$$\lim_{n\to\infty} \frac{S(n)}{n} = \lim_{n\to\infty} M(x_n) = 3/\pi^2 . \tag{4.7}$$

In this case, the stationarity condition (2.6) in the broad sense is satisfied.

Taking into account that the random variables of the corresponding sequence are bounded $|x_n| \leq 2$, then the stationarity condition (2.4) in the broad sense is satisfied.

The asymptotic independence of the random sequence $x_n$, ie, the condition (2.3) is satisfied for an arithmetic function for which the presence of an even or odd number of prime divisors is the defining property.

Thus, the whole stationarity condition in the broad sense is satisfied for the corresponding $S(n)$ sequence of random discrete quantities $x_n (n=1,2,...)$.

Based on (4.5) and (4.6), we find the variance of the random variable $x_n$:

$$D(x_n) = M(x_n^2) - M^2(x_n) = 15/\pi^2 - 9/\pi^4 + o(1) \neq 0. \tag{4.8}$$

Consequently, in order to verify sufficient conditions for the limit theorem it remains to verify the strong mixing property with $\sum_{n=1}^{\infty} \alpha(n) < \infty$ for sequences $x_n (n=1,2,...)$ corresponding to the indicated arithmetic functions $S(n)$ .

We take the event $A$ that the number $k$ has an even number of prime divisors of the first degree, i.e. $x_k = 2$. Let's take the event $B$ that the number $n+k$ is not free squares, i.e. $x_{n+k} = 0$.

Then $P(x_k = 2) = 3/\pi^2 + o(1)$, $P(x_{n+k} = 0) = 1 - 6/\pi^2 + o(1)$. In this case the product of probabilities is equal to: $P(A)P(B) = [3/\pi^2 + o(1)][1 - 6/\pi^2 + o(1)]$. Then the value $P(AB)$ is: $P(A)P(B/A) = [3/\pi^2 + o(1)][1 - 6/\pi^2 + o(1)]$. The last equality is explained by the independence of the event $B$ from the event $A$. . Therefore, in this case, the condition of strong mixing is satisfied:

$$\alpha(n) = \sup |P(AB) - P(A)P(B)| = 0. \tag{4.9}$$

Based on (4.9): $\sum_{n=1}^{\infty}\alpha(n)=0<\infty$, therefore, all the conditions of the above limit theorem hold for a sequence of random variables corresponding to an arithmetic function $S(n)$, and the limiting distribution for this function is normal.

It was shown that the Mertens function satisfies all stationary conditions in the broad sense. In addition, random variables $x_n$ corresponding to the Mertens function are bounded - $|x_n|\leq 1$. A sufficient condition of the limit theorem for variance is also satisfied: $D(x_n)=6/\pi^2+o(1)\neq 0$.

We prove the fulfillment of a sufficient condition for the limit theorem on strict mixing for an arithmetic function $M(n)$ similarly, as for an arithmetic function $S(n)$. Consequently, the limiting distribution for the Mertens function $M(n)$ is also normal. This is consistent with the QEIS schedule A002321.

## 5. AN ERGOTIC THEOREM FOR ARITHMETIC FUNCTIONS

We have considered the summation arithmetic functions in Section 2, which are the sum of a stationary sequence of random variables $x_i (i=1,2,...)$ in the broad sense.

Now we consider $H^2 = H^2(\Omega, F, P)$ - the space of (complex-valued) random variables $x = a + bi, a, b \in R$, with $D(x) < \infty$, where $|x|^2 = a^2 + b^2$.

If $x, y \in H^2$, then we assume that the scalar product is equal to:

$$(x,y)=M(x\bar{y}),\tag{5.1}$$

where $\bar{y}=a-bi$, is the complex conjugate of $x=a+bi$ and the norm:

$$\|x\|=(x,x)^{1/2}.\tag{5.2}$$

The space $H^2$ is a complex Hilbert space (of random variables considered on a probability space $(\Omega, F, P)$).

Having in mind (5.1), (5.2) it follows that if $M(x)=M(y)=0$, then:

$$\mathrm{cov}(x,y)=(x,y).\tag{5.3}$$

Note that the independent random variables x, y in (5.3) correspond to orthogonal random variables (the scalar product, which is equal to zero - $cov(x, y) = 0$) in a particular case.

Let's consider a sequence of random variables of the form:

$$y_n = \sum_{k=0}^{\infty} a_k \xi_{n-k}, \qquad (5.4)$$

where $a_k$ are complex numbers for which $-\sum_{k=0}^{\infty} |a_k|^2 < \infty$, and it is performed $M(\xi_{n-k}) = 0, D(\xi_{n-k}) = 1, M(\xi_i, \xi_j) = 0$ for random variables $\xi_{n-k}$, if $i \neq j$.

Thus, random variables $\xi_i$ form an orthonormal basis along which a random variable $y_n$ is decomposed.

Any sequence of independent random variables $x_n$ can be represented in the form:

$$x_n = m_n + y_n, \qquad (5.5)$$

where $m_n$ is the mathematical expectation $x_n$, and $y_n$ is determined by (5.4).

Therefore (5.5) is a Wald expansion for an independent random sequence $x_n$: $m_n$ is a singular sequence, and $y_n$ is a regular sequence.

Now consider another sequence of complex random variables:

$$x_n = \sum_{k=-\infty}^{\infty} z_k e^{i\lambda_k n}, \qquad (5.6)$$

where $z_k$ are the orthogonal random variables ($M(z_i \bar{z}_j) = 0, i \neq j$) with zero means and variances $D(z_k) = \sigma_k^2 > 0, -\pi \leq \lambda_k \leq \pi, \lambda_i \neq \lambda_j, i \neq j$. We can understand $z_k$ (as previously mentioned) independent random variables in the particular case.

If we put:

$$Z(\lambda) = \sum_{\{k, \lambda_k \leq \lambda\}} z_k, \qquad (5.7)$$

then (5.6) based on (5.7) can be written in the form:

$$x_n = \sum_{k=-\infty}^{\infty} e^{i\lambda_k n} \Delta Z(\lambda_k), \qquad (5.8)$$

where $\Delta Z(\lambda_k) = Z(\lambda_k) - Z(\lambda_k -) = z_k$.

The right-hand side in (5.8) is the integral sum for the stochastic integral. Therefore, having in mind (5.8), the sequence (5.6) can be written in the form of the following stochastic integral:

$$x_n = \int_{-\pi}^{\pi} e^{i\lambda n} dZ(\lambda), \qquad (5.9)$$

where $Z(\lambda)$ is the process with orthogonal increments.

We know an expression analogous to (5.9) for the covariance function (Herglotz's theorem):

$$R(n) = \operatorname{cov}(x_n x_o) = M(x_n x_0) = (x_n, x_0) = \int_{-\pi}^{+\pi} e^{i\lambda n} dF(\lambda), \qquad (5.10)$$

where $F(\lambda)$ is the spectral function. The process with orthogonal increments $Z(\lambda)$ satisfies the relation in this case:

$$M \, |dZ(\lambda)|^2 = dF(\lambda). \qquad (5.11)$$

The ergodic theorem (in the mean square sense) is proved in [8].

Let the sequence $x_n$ is stationary, in a broad sense, with $M(x_n) = 0$ and covariance function $R(n)$ with spectral decomposition (5.10). Then:

$$\frac{1}{n} \sum_{k=0}^{n-1} x_k \xrightarrow{L^2} Z(\{0\}), \qquad (5.12)$$

$$\frac{1}{n} \sum_{k=0}^{n-1} R(k) \to F(\{0\}). \qquad (5.13)$$

Corollary 1. If the spectral function $F(\lambda)$ is continuous at zero, i.e. $F(\{0\}) = 0$, then, almost certainly $Z(\{0\}) = 0$, and (5.12), (5.13) will be written in the form:

$$\frac{1}{n}\sum_{k=0}^{n-1} x_k \xrightarrow{L^2} 0, \qquad (5.14)$$

$$\frac{1}{n}\sum_{k=0}^{n-1} R(k) \to 0. \qquad (5.15)$$

Corollary 2. If the stationary sequence in the broad sense $x_n$ has a mathematical expectation $M(x_n) = m \neq 0$ and the function $F(\lambda)$ is continuous at zero, then (5.12), (5.13) can be written in the form:

$$\frac{1}{n}\sum_{k=0}^{n-1} x_k \xrightarrow{L^2} m, \qquad (5.16)$$

$$\frac{1}{n}\sum_{k=0}^{n-1} R(k) \to 0. \qquad (5.17)$$

It is proved [5] that if a stationary in the wide sense sequence of random variables is regular, then the spectral function $F(\lambda)$ is absolutely continuous. Therefore, based on Corollary 1, the limiting relations (5.14) and (5.15) hold in this case.

Based on (5.5) the sequence of independent random variables $x_n$ is regular if $M(x_n) = 0$. Therefore, the relations (5.14) and (5.15) hold for the given sequence. If the sequence of independent random variables has a mathematical expectation $M(x_n) = m \neq 0$, then, based on Corollary 2, the limiting relations (5.16) and (5.17 hold for the sequence.

I recall that we are interested in the behavior of arithmetic functions for $n \to \infty$, i.e. we investigate the behavior of the summation arithmetic functions $S(n)$, which are the sum of independent random variables $x_n$ for $n \to \infty$.

Thus, having in mind (5.16), we can write the following limiting relation for the summation arithmetic functions, which are the sum of a stationary, in a broad sense, sequence of independent random variables $x_n$ for $n \to \infty$ with a mathematical expectation $m = \lim_{n \to \infty} M(x_n)$:

$$\frac{S(n)}{n} \xrightarrow{L^2} m. \qquad (5.18)$$

Expression (5.18) will be called an ergodic theorem for summation arithmetic functions.

It follows from (5.18):

$$S(n) - nm \xrightarrow{L_2} 0. \tag{5.19}$$

Based on (5.19), the standard deviation of the quantity $S(n) - nm$ is:

$$\sigma(S(n) - nm) = O(n^{1/2}). \tag{5.20}$$

# 6 ASYMPTOTICS OF THE GROWTH OF THE DEVIATION OF ARITHMETIC FUNCTIONS FROM THE MEAN VALUE

Following [3], on the basis of Chebyshev's inequality, we can write the relation for any real arithmetic function $f(n)$ on the probabilistic space considered above:

$$\lim_{n \to \infty} P(|f(n) - M_f(n)| \leq \Psi(n)\sqrt{D_f(n)}) = 1, \tag{6.1}$$

where $M_f(n), D_f(n)$, respectively, the mean value and variance $f(n)$, and $\Psi(n)$ is an arbitrary slowly growing function. Expression (6.1) is some analogue of the law of large numbers for an arithmetic function $f(n)$.

Having in mind (6.1), almost everywhere, the following relation holds:

$$|f(n) - M_f(n)| \leq \Psi(n)\sqrt{D_f(n)}, \tag{6.2}$$

where as a function we can take any slowly increasing function $\Psi(n)$.

It follows from (6.2), that a definition $D_f(n)$ is very important to determine the deviation from the mean value of an arithmetic function..

The estimate for the variance of an arithmetic function $S(n)$ that is the sum of independent or asymptotically independent bounded random variables $x_i$ is:

$$D_S(n) = D(\sum_{i=1}^{n} x_i) = \sum_{i=1}^{n} D(x_i) \leq n \sup_i \{D(x_i)\}. \tag{6.3}$$

Based on (6.2), (6.3), we obtain an estimate of the deviation of the arithmetic function from its mean value for the case:

$$|S(n) - M_S(n)| \leq \Psi(n)\sqrt{n \sup_i \{D(x_i)\}}. \tag{6.4}$$

As was shown above, the relation for an arithmetic function $Q(n)$ that is the sum of independent random variables are fulfilled - $D_f(n) = np(1-p)$. Therefore, the estimate for this case is:

$$D_f(n) \leq n/4. \tag{6.5}$$

Based on (6.4) and (6.5), almost everywhere an estimate is performed in this case:

$$|Q(n) - M_Q(n)| \leq 0,5 n^{1/2} \Psi(n) \tag{6.6}$$

There is an estimate [5] of the variance value for the summation arithmetic functions $D(S_n)$, which are the sum of dependent random variables, if $\lim_{n \to \infty} D(S_n) = \infty$, when the stationarity conditions in the broad sense and the strong mixing are fulfilled for a sequence of random variables $x_n (n = 1, 2, ...)$:

$$D(S_n) = n h(n), \tag{6.7}$$

where $h(n)$ is a slowly growing function of the natural argument, i.e. in this case $D(S_n)$ grows almost linearly.

Having in mind (6.2) and (6.7), the following relation holds, almost everywhere, in this case:

$$|S(n) - M_S(n)| \leq h_1(n) n^{1/2}, \tag{6.8}$$

where $h_1(n)$ is a slowly growing function of the natural argument.

The following estimate holds for the function $h_1(n)$:

$$h_1(n) \leq n^\xi, \tag{6.9}$$

where $\xi$ is a small real positive number.

Therefore, based on (6.8), (6.9), the relation for a given deviation is valid almost everywhere, if the stationarity conditions in the broad sense and the strong mixing are satisfied for the corresponding sequence of random variables $x_n (n = 1, 2, ...)$:

$$|S(n) - M_S(n)| \leq n^{1/2+\xi}. \tag{6.10}$$

An estimate of the deviation of the arithmetic function from the mean value (6.10) is satisfied for all the non-derivatized arithmetic functions considered above.

For example, for the Mertens function $M(n)$, the average value is 0, therefore, on the basis of (6.8), almost everywhere performed:

$$|M(n)| \leq n^{1/2+\xi}, \qquad (6.11)$$

which corresponds to (5.20) and the equivalent formulation of the Riemann hypothesis.

## 7. CONCLUSION AND SUGGESTIONS FOR FURTHER WORK

The next article will continue to study the asymptotic behavior of arithmetic functions.

## 8. ACKNOWLEDGEMENTS

Thanks to everyone who has contributed to the discussion of this paper. I am grateful to everyone who expressed their suggestions and comments in the course of this work.

# References

1. Vinogradov I.M. Osnovy teorii chisel [Fundamentals of number theory]. M.: Nauka [Moscow: Publishing House «Science»]. 1981.

2. Karatsuba A.L. Osnovy analiticheskoy teorii chisel [Fundamentals of analytic number theory]. M.: Nauka.Glavnaya redaktsiya fi ziko-matematicheskoy literatury [Moscow: Publishing House «Science». Home edition of physical and mathematical literature]. 1983. 240 p.

3. Kubilyus Y. Veroyatnostnye metody v teorii chisel [Probabilistic methods in number theory]. Vilnius, 1962. 220 p.

4. Hausman M., Shapiro H. Distribution functions of huperadditive arithmetic functions. 1987. Vol. 40. March. Pp. 221…227.

5. Ibragimov I.A., Linnik Yu.V. Nezavisimye i statsionarno svyazannye velichiny [Independent and permanently connected quantities]. M.: Nauka [Moscow: Publishing House «Science»]. 1965. 524 p.

6. Bukhshtab A.A. Teoriya chisel [Number theory]. M.: Prosveshchenie [Moscow: Publishing House «Education»]. 1966. 384 p.

7. Volfson V.L. Investigation of the asymptotic behavior of the Mertens function, arXiv preprint https://arxiv.org/abs/1712.04674 (2017)

8. Shiryaev A.N. Veroyatnost-2 [Probability-2]. Moscow: MTsNMO Publishing House, 2007. 574 p.